\documentstyle{article}

\newcommand{\R}{{\bf R}}
\newcommand{\os}{{\overline{\Cc}_X}}
\newcommand{\op}{{\overline{\Pp}}}
\newcommand{\tsi}{{\widetilde\si}}
\newcommand{\ttau}{{\widetilde\tau}}
\newcommand{\TX}{{\widetilde X}}
\newcommand{\Tilde}{\widetilde}
\newcommand{\id}{{\rm id}}
\newcommand{\Gr}{{\rm Gr\,}}

\newcommand{\Q}{{\bf Q}}

\newcommand{\Nn}{{\cal N}}
\newcommand{\Cc}{{\cal C}}
\newcommand{\Z}{{\bf Z}}
\newcommand{\C}{{\bf C}}

\newcommand{\Ee}{{\cal E}}
\newcommand{\al}{{\alpha}}

\newcommand{\be}{{\beta}}

\newcommand{\om}{{\omega}}
\newcommand{\eps}{{\varepsilon}}
\newcommand{\de}{{\delta}}

\newcommand{\ga}{{\gamma}}

\newcommand{\ka}{{\kappa}}
\newcommand{\la}{{\lambda}}
\newcommand{\si}{{\sigma}}
\newcommand{\Pp}{{\cal P}}
\newcommand{\PD}{{\rm PD}}

\newcommand{\Emb}{{\it Emb}}

\newcommand{\Si}{{\Sigma}}

\newcommand{\SS}{{\smallskip}}
\newcommand{\MS}{{\medskip}}

\newcommand{\NI}{{\noindent}}
\newcommand{\proof}[1]{\noindent{\bf Proof#1:\  }}

\newcommand{\QED}{\hfill$\Box$\medskip}

\newtheorem{theorem}{Theorem}[section]
\newtheorem{cor}[theorem]{Corollary}

\newtheorem{remark}[theorem]{Remark}
\newtheorem{lemma}[theorem]{Lemma}

\newtheorem{prop}[theorem]{Proposition}

\newcommand{\at}{{@}}

\title{From  symplectic deformation to  isotopy}
\author{Dusa McDuff\thanks{Partially supported by
NSF grant DMS 9401443.} \\ State University of New York at Stony Brook \\
(dusa\at math.sunysb.edu)}

\date{May 15, 1996, revised Aug 4, 1997}

\begin{document}

\maketitle

\begin{abstract}  Let  $X$ be an oriented
 $4$-manifold which does
not have simple SW-type, for example a blow-up of a rational or ruled
surface.  We show that any two cohomologous and deformation equivalent
symplectic forms on $X$  are isotopic.  This
implies that blow-ups of these manifolds are unique, thus extending
work of Biran.  We also establish uniqueness of structure for certain fibered
$4$-manifolds. \end{abstract}

\section{Introduction}

 Two symplectic
forms $\om_0, \om_1$ on $X$ are said to be {\bf deformation equivalent}
if they may be joined by a family of  symplectic
forms,    and are called {\bf isotopic} if this family may be chosen so that
its elements all  lie in the same cohomology class.   Moser showed that every
isotopy on a compact manifold has the form $\phi_t^*(\om_0)$, where $\phi_t:
X\to X$ is a family of diffeomorphisms starting at the identity.  Examples are
known in dimensions $6$ and above of cohomologous symplectic forms that are
deformation equivalent but not isotopic: see~\cite[Example~7.23]{INTRO}.  No
such examples are known in dimension $4$, and it is a possibility that the two
notions are the same in this case.   (Note that the Gromov invariants are
deformation invariants -- in fact, Taubes's recent work shows that they are
smooth invariants -- and so they cannot help with this question.)

In~\cite{LM} Lalonde--McDuff developed ideas from~\cite{LAL,RR} into
  an \lq\lq inflation" procedure
that converts a deformation into an isotopy, and
 applied it to
establish the  uniqueness of symplectic structures on ruled surfaces. The
present note extends the range of
this procedure and describes various applications of it.

Here is the basic lemma:
its proof is sketched in~\S2 below.  
We will denote the Poincar\'e dual of a class $A$ by
$\PD(A)$, and will express the result in terms of the Gromov invariant
 $\Gr_0(A)$ that
counts {\it connected} $J$-holomorphic representatives 
of the class $A$.  (This is the invariant used by Ruan in~\cite{RU}. To
 a first approximation, it is the same as the 
invariant used by Taubes.)

\begin{lemma}[Inflation Lemma]\label{infl}  Let $A$ be a  class in $
H_2(X,\Z)$, with  nonnegative
self-intersection number and nonzero Gromov invariant $\Gr_0(A)$.
Moreover, if $A^2 = 0$
assume that $A$ is a primitive class.  
Then, given any family $\om_t, 0\le t\le 1,$ of
symplectic forms on $X$ with $\om_0 = \om$, 
there is a family $\rho_t$ of closed
forms on $X$ in class $\PD(A)$ such that the family  $$ 
\om_t + \ka(t)\rho_t,\quad 0 \le t\le 1,
$$
is symplectic whenever $\ka(t) \ge 0$.
\end{lemma}

To apply this lemma, we need to understand the Gromov
invariants of $X$.
 Recall that a  symplectic $4$-manifold $X$ is said to have
{\bf simple SW-type} or just {\bf simple type} if its only nonzero Gromov
invariants occur in classes $A\in H_2(X)$ for which
$$
k(A) = - K\cdot A + A^2 = 0.
$$
Taubes showed in~\cite{TAU0,TAU00} that all symplectic $4$-manifolds
with $b_2^+ > 1$ have simple type.  It follows easily from the wall crossing
formula of Li--Liu~\cite{LL1,LL2} that if $b_2^+ = 1$ then $X$ has simple
type only if $b_1 \ne 0$ and  all products of elements in $H^1(X)$ vanish.
Moreover, Liu showed in~\cite{LIU} that any minimal symplectic manifold
with $K^2 < 0$ is ruled.  Therefore the symplectic $4$-manifolds with
nonsimple type are blow-ups of  \MS

\NI
(i)  rational and ruled manifolds; \SS

\NI
(ii)
 manifolds (such as the Barlow or Enriques surface) with $b_1 = 0, b_2^+ = 1$;
\SS

and
\SS

\NI
(iii) manifolds with $b_1 = 2$ and $(H^1(X))^2 \ne 0$.  Examples 
with $K = 0$ are
hyperelliptic surfaces and some (but not all) of the  non-K\"ahler
$T^2$-bundles over $T^2$.  There are
also K\"ahler examples with $K\ne 0$,	 for
example  quotients of the form  $T^2\times \Si/ G$ where $\Si$ 
has genus $> 1$ and $G$ is a suitable finite group.\footnote
{
Le Hong Van pointed out that such examples exist: they were omitted
from the survey article~\cite{MS}.}

\SS
  We will show below that
manifolds which do not have simple type have enough nonzero Gromov
invariants for the following result to hold.

\begin{theorem}\label{1} Let $(X, \om)$ be a symplectic $4$-manifold which
does not have simple type.  Then any deformation between two
cohomologous symplectic forms on $X$ may be homotoped through
deformations with fixed endpoints  to  an isotopy.
\end{theorem}

\begin{cor}  For any $k > 0$ and  any $X$ not
of simple type, there is at most one way of blowing up $k$
points to specified sizes.
\end{cor}
\proof{}  It is easy to see that any two blow-ups
of a symplectic manifold are deformation equivalent. (For more details of this
step see the proof of Corollary~\ref{c2} below.)  The cohomology class of the
blow-up is determined by the \lq\lq size" of the blown-up points.   Hence the
above theorem implies that cohomologous blow-ups are isotopic.\QED

When $k = 1$ this corollary was first proved by McDuff~\cite{BL,INTRO} in
the case of $\C P^2$ and ruled surfaces over $S^2$,
and by Lalonde~\cite{LAL} for (certain) irrational ruled surfaces.  The
methods used were very geometric.  Recently, Biran~\cite{BIR}  established
uniqueness  for many blow-ups of $\C P^2$.  Independently,
Gatien~\cite{GAT}
realised that  arguments similar to the ones in this note should lead to the
uniqueness of blow-ups  for  ruled manifolds.

The next result was proved by
Lalonde--McDuff in~\cite{LM}.

\begin{cor}\label{rruniq}  Let $X$ be oriented diffeomorphic to a 
minimal rational or
ruled surface, and let $a\in H^2(X)$.  Then there is a symplectic 
form on $X$ in
the class $a$ which is compatible with the given orientation if and only if
\,\,$a^2 > 0$.  Moreover, any two symplectic forms in class \,\,$a$ are
diffeomorphic. \end{cor}
\proof{}
The first result is obvious,  and the second holds  because, by~\cite{RR},
given
any two symplectic forms $\om_0, \om_1$ on $X$ there is a diffeomorphism
$\phi$ of $X$ such that $\phi^*(\om_0)$ is deformation equivalent to
$\om_1$.\QED

\subsubsection*{Embedding balls and blowing up}

Using  the correspondence between embedding balls and
 blowing up that was first described in~\cite{BL}, we can
deduce the connectedness of the space
$$
\Emb(\coprod_{i=1}^k\,B(\la_i), X)
$$
of
symplectic embeddings of the disjoint union of the
closed $4$-balls $B(\la_i), i =
1,\dots, k,$ of radius $\la_i$ into $X$, with the $C^1$-topology.

\begin{cor}\label{c2} If $X$ has nonsimple type,
 $\Emb(\coprod_{i=1}^k\,B(\la_i),X)$ is path-connected.
\end{cor}

The proof is sketched later.  Note that a similar result holds for manifolds of
the form $X - Z$ where $X$ has nonsimple type as above and
 $Z$ is any symplectic
submanifold.  In particular, it holds for the open unit ball, which is just
$\C P^2 - \C P^1$.   

\begin{remark}\rm
Here we have applied the inflation procedure to discuss the 
uniqueness problem for blow-ups.  Closely related is the packing problem:  
what constraints are there on the radii of balls that embed symplectically
and disjointly in $X$?  The paper \cite{MCDPOL} gives a 
complete solution of this
problem for $k\le 9$ balls in $\C P^2$, but left open the 
question of whether there is
a full filling (i.e. a filling that uses up 
all the volume) of $\C P^2$ by $k$ equal balls
for $k\ge 10$.  Using the methods of the 
present paper, Biran has shown that such a
filling does indeed exist: see~\cite{BIR2}.
\end{remark}

\MS

\subsubsection*{The symplectic cone}

Our methods also give some information for  general symplectic $4$-manifolds.
As Biran points out in his discussion of the packing problem,  one can use
Lemma~\ref{infl} to get information on the symplectic cone
$$
\Cc_X = \{a\in H^2(X;\R): a\mbox{ has a symplectic representative }\}.
$$
(Here we consider only those symplectic 
forms that are compatible with the given
orientation of $X$.)
Because, by Taubes, the set of classes
with nontrivial Gromov invariants is an 
invariant of the smooth structure of $X$,
it is clear that
 $$
\Cc_X\subset\{ a : a(A) > 0 \mbox{ for all }
  A\ne 0\mbox{ such that }\Gr(A)\ne 0\}.
$$

\begin{prop} 
(i)  If there is some
symplectic form $\om$ on $X$ for which the class $A$
 is represented by a  symplectic submanifold with all
components of nonnegative self-intersection, 
then $\PD(A)$ is in the closure $\os$ of
the symplectic cone $\Cc_X$.  Moreover, if
 $X$ has simple type and $\Gr(A)\ne 0$ for
some class  $A\ne K, 0$, then $\PD(A)\notin \Cc_X$.  
\SS

\NI(ii)
Suppose that  $X$ is minimal with $b_2^+ > 1$ and that $K$ is nontorsion.  Then
$\PD(K)\in\os$, while $K\notin\Cc_X$ if there is any nonzero class $A$
with $A^2 = 0$ and $\Gr(A) \ne 0$.  \end{prop}
\proof{}  The first statement in (i) follows immediately from the Inflation Lemma, 
because  the class $\frac 1{\ka} [\om + \ka \rho]$
converges to $[\rho]$ as $\ka \to \infty$.\footnote{
We do not need
any hypothesis here about nonmultiply covered toral components since 
we are applying this lemma  using a single submanifold $Z$ rather than a family
$Z_t$.}
If $\Gr(A)\ne 0$ then the fact that $X$ has simple type immediately implies
that $A\cdot (K - A) = 0$.  Moreover, Taubes showed that $\Gr(A) = \pm \Gr(K-A)$. 
Therefore, if both $A$ and $K-A$ are nonzero, $\PD(A)\notin\Cc_X$.  
(ii) follows immediately from Taubes' results since the hypotheses
imply that $K\cdot A = 0$.
\QED

It would be interesting to know whether $K\in \Cc_X$ in the case when 
$K$ has a representative with all components of positive self-
intersection.

\subsubsection*{Deformations on arbitrary $X$}

 Let us write $\Ee$
for the subset of $H_2(X,\Z)$ consisting of all classes that can be represented
by a symplectically embedded $2$-spheres of self-intersection $-1$.    If $X$
is
not the blow-up of a rational  or ruled manifold, it was shown in
McDuff~\cite{IMM} that $\Ee$ is finite and consists of mutually orthogonal
classes, ie $E\cdot E' = 0$ for $E\ne E' \in \Ee$.  In particular, $X$ is the
symplectic blow-up of a unique minimal manifold $Y$ which is obtained by
blowing down a set of representatives of the classes in $\Ee$.
Now, the Seiberg--Witten blow-up formulas imply that the
 Gromov invariants of $X$ are determined by those of $Y$ by the following
rule:\footnote
{
Note that the theory of $J$-holomorphic curves implies only
that  $\Ee$ and hence $Y$ depend  on the deformation class of $\om$, while
Taubes--Seiberg--Witten theory shows that  when $b_2^+ >
1$ they are  smooth invariants of $X$.}
there is a  natural decomposition
 $$
H_2(X;\Z) = H_2(Y;\Z) \oplus \sum_{E\in\Ee} \Z E,
$$
and for each
class $B\in H_2(Y;\Z)$
\begin{eqnarray*}
\Gr_X(B +\sum_{E\in\Ee}  \eps_E E) & = & \Gr_Y(B)
\quad\mbox{ if } \;\eps_E
= 0\mbox{ or }\;1 \;\;\mbox{ for all $E$},\\
& = & 0\quad\mbox{otherwise}.
\end{eqnarray*}

Taubes showed in~\cite{TAU00} that each $J$-holomorphic representative
of the class $B $ consists  of a finite number  of
disjoint components with nonnegative self-intersection, and that the class
$B +\sum \eps_E E$ is represented by adjoining  the appropriate exceptional
curves.  Moreover, these components are embedded, except possibly if they
represent non-primitive classes $T$ with $T^2 = 0$, in which case they may be
multiply-covered  tori of self-intersection zero.  Since these multiply-covered
tori do not always
persist under deformation (see~\cite{TAUTOR}), we cannot use them.
Therefore, as before, we use the invariant
 $\Gr_0(A)$, which is defined to be the
algebraic number of {\it connected} $J$-holomorphic representatives of the
class $A$ if either $A^2 > 0$ or $A$ is primitive. 
 We then define
$$
V = V(X) \subset
H^2(X,\R)
$$
to be the convex cone spanned by the Poincar\'e duals of those
classes $A$ such that
\SS

\NI (i)  $A\cdot E \ge 0$ for all $E\in \Ee$;

\NI (ii) either  $A$ is primitive or $A^2 > 0$; and

\NI (iii)  $\Gr_0(A)\ne 0$.
\SS

\begin{remark}\rm
It is not hard to see  (and details may be found in~\cite{MONT}) that
 if $A$ has no $J$-holomorphic 
representatives that include components that are multiply-covered tori then
$$
\Gr(A) = \sum_{A_1+\dots + A_k = A} \prod_i \Gr_0(A_i),
$$
where the sum is taken over all sets $\{A_1,\dots, A_k\}$ such that
 $A = A_1 + \dots A_k$.
Further, if all regular
multiply-covered $J$-holomorphic tori in $X$  lie in classes which are
multiples of primitive classes with nonzero Gromov invariants, then $V$ may
be defined by replacing the conditions (ii) and (iii) above by
\SS

\NI (iv)  $\Gr(A) \ne 0$.
\end{remark}

Given these definitions, the inflation procedure
 applies immediately to show:

\begin{theorem}\label{2} Let $X$ be any symplectic $4$-manifold.  If $\om_t,
t\in [0,1],$ is a deformation with $[\om_0] = [\om_1]$ such that
$$
[\om_0] \in  \R^+[\om_t] + V
$$
 for all $t$, then $\om_t$
can be homotoped (rel endpoints) to an isotopy.
\end{theorem}

The above theorem does not give any interesting information about
blow-ups.   For $V$ lies in the annihilator of $\Ee$ so that we cannot
now change the size of the exceptional curves. 
However, as we now show, one can sometimes  combine this result
with other geometric information about the symplectic 
structure to get something new. 

\subsubsection*{Symplectic fibrations}

Let $\pi: X\to B$ be a fibration with compact oriented total space $X$ and
oriented base $B$. A symplectic form $\om$ on $X$ is said to be
$\pi$-compatible 
 if all the fibers of $\pi$ are symplectic
submanifolds of $X$ and if the orientations that $\om$ defines on $X$
and the fibers equal the given ones.  We will
see in \S3 that when the base and fiber have dimension $2$ all such
forms are deformation equivalent.  Our methods allow us to change
this deformation into an isotopy when $X$ has nonsimple type (for
example, if the base or fiber is a sphere or if $X$ is a hyperelliptic
surface) and also in some other cases.  Here is a sample result that
applies when $X = F\times B$. We write $\si_F, \si_B$
for the pullback to $X$ of area forms on $F$ and $ B$ 
that have total area $1$ and $g_F, g_B$ for the genus of $F$ and $B$.

\begin{prop}\label{fib}  Suppose that $\om$ is a $\pi$-compatible form on 
the product $X = F\times B$ where $g_F > 1$.  If $[\om]
= [\si_F + \la\si_B]$, where 
$\la (g_F - 1) > g_B-1$, then  $\om$ 
is isotopic to the split form $\si_F + \la\si_B$.\end{prop}

\NI
The proof is given in \S3.
\SS

 \NI
{\it Acknowledgement}  I wish to thank Paul Biran and Francois Lalonde for 
useful comments on an earlier draft of this note.

\section{Inflation and manifolds of nonsimple type}

We begin by sketching the proof of the inflation lemma.

\subsection*{Proof of Lemma~\ref{infl}}

Let $\om_t$ be any
deformation with $[\om_0] = [\om_1]$, and choose a generic  family $J_t$ of
$\om_t$-tame almost complex structures. Lalonde and McDuff show
in~\cite{LM} that if $A$ with $\Gr(A)\ne 0$ 
may be represented by a symplectically embedded
submanifold of positive self-intersection then, for some smooth map  $t\mapsto
\mu(t)$ with $\mu(0) = 0, \mu(1) = 1$,  there is a smooth family $Z_t$ of
$J_{\mu(t)}$-holomorphic embedded submanifolds of $X$ in class $A$.
The same argument holds when $A^2 = 0$ provided that $A$ is primitive.
\footnote
{
If $A^2= 0$ and $A$ is not primitive, this may not hold because of problems
in counting multiply-covered tori: see Taubes~\cite{TAUTOR}.}
  Without loss of generality we
may suppress the reparametrization, and 
suppose that $\mu (t) = t$.  It is shown in
Lemma~3.7 of~\cite{RULED} that there is a smooth family $\rho_t$ of closed
$2$-forms on $X$ that represent the class $\PD(A)$ and are such that
$$
\om_t+ \ka(t) \rho_t
$$
is symplectic for all constants $\ka(t) \ge 0$ and all $t \in
[0,1]$.  For completeness, we sketch this construction for 
fixed $t$, noting that the
construction can be made to vary smoothly with $t$. 

The form $\rho$ 
is supported near $Z$ and
represents the Thom class of the normal bundle to $Z$.
  If this normal bundle is
trivial, a
neighborhood of $Z$ is symplectically equivalent to a product $Z\times
D^2$ and the existence of $\rho$ is obvious.   In general, let  $k = Z\cdot Z$
and choose a connection $\ga$ on the normal circle bundle
$\pi:Y \to Z$ such that
$d\ga  = -\pi^*(f\om_Z)$, where $\om_Z = \om|_Z$ and the function $f\ge 0$ has 
 appropriate integral over $Z$.  Then a
neighborhood of $Z$ can be symplectically identified with a neighborhood
of the zero section in the associated complex line bundle,  equipped
with the symplectic form
$$
\tau = \pi^*(\om_Z) + d(\pi{r^2}\ga),
$$
where $r$ is the radial distance function.  
Hence one can take $\rho$ to be given
by the formula $-d(g(r)\ga)$, where $g(r)$ is a nonnegative
function with support in a small
interval $[0,\eps]$ that equals $  \pi {r^2} - 1$
 near $r= 0$.
Note that if $Z$ had negative
self-intersection, one would have to take $f \le 0$, thus making
the integral
of $\rho$  over $Z$ negative.  Hence in this case $\om_t+ \ka(t)
\rho $ would cease to be symplectic for large $\ka(t)$.
\QED

\begin{remark}\label{many}\rm {\bf (i)}
Once the family $Z_t$ has been found we can alter $J_t$, keeping
the $Z_t$ $J_t$-holomorphic, so that $\rho_t$ is semi-positive, i.e.
$$
\rho_t(x, J_t x) \ge 0, 
$$
for all tangent vectors $x$.  (One can achieve this by making  both the
fibers of the normal bundle and their  orthogonal complements with respect 
to the closed form $d(r^2\ga)$
 $J_t$-invariant.)  It follows that all the forms $\om_t + \ka(t)\rho_t$
tame $J_t$.
\SS

\NI
(ii)
We will have to apply this process repeatedly, along families of submanifolds
$Z_{1t},\dots, Z_{nt}$ which intersect.  The procedure as outlined above
depends on the order chosen for the $Z_{it}$.  However,
 with a little more care one can make the 
processes commute.  (Alternatively, we
can perform all inflations simultaneously.)  
Since we are dealing with a $1$-parameter deformation, we may suppose that 
for each $t$ the
manifolds  $Z_{jt}, j = 1,\dots, n,$ meet transversally in pairs.  
If  these intersections are all $\om_t$-orthogonal,
 then it is not hard to see that we
can change   the forms $\rho_{it}$ so that
$$
\om_t + \sum_i \ka_i(t) \rho_{it}
$$
is symplectic for all $\ka_i(t) \ge 0$.  
Indeed, to do this we just need to be careful
near intersection points and here the local model is a product
$U_i\times U_j$  with a product form, where $U_i\subset Z_{it}$.  
Therefore, if we choose the functions $f_i$ above so that they
vanish on $U_i$, the  normal form $\tau$ for $\om$ will
respect this splitting, as will the forms $\rho_{it}, \rho_{jt}$.
xs

To arrange that the intersections are orthogonal, one needs to perturb
 the families $Z_{it}$. Remark
that two kinds of perturbation are needed here.  By
 positivity of intersections, every
intersection of $Z_{it}$ with $Z_{jt}$ counts positively and so, 
by a  $C^1$-small
perturbation we can arrange that these intersections are  transverse for all
$i,j,t$.  (Observe that singularities can always be avoided for 
$1$-parameter families,
but possibly not for higher dimensional families.)  We then need
 a large perturbation
to make the intersections symplectically orthogonal:  this can be 
done by the methods
of~\cite{RULED}. Note that after these perturbations the 
resulting manifolds need no
longer be $J_t$-holomorphic, though they will be holomorphic
 for some other family
$J_t'$. \SS

\NI
{\bf (iii)}
Another way of viewing this inflation process is as a form of the Gompf sum.
To inflate along a connected submanifold $Z$ one simply identifies $X$ with a
sum of the form\footnote
{This identification of $X$ with $X\#_{Z = Z_-}W$ is a version
of the  \lq\lq thinning" process of~\cite{MSYM}.}
$$
X\#_{Z = Z_-}W,
$$
where $W$ is a ruled surface over $Z$ which has a section $Z_-$ with
self-intersection equal to $-Z^2$, and then inflates $X$ by increasing the size
of the fiber of $W$.   When $Z$ has different
components that intersect $\om$-orthogonally one needs to replace the Gompf
sum by the kind of plumbing process used by McDuff-Symington
in~\cite{MSYM}. \end{remark}

\subsection*{Gromov invariants of manifolds of nonsimple type}

The following lemma contains all the information we need on Gromov
invariants.  We will denote the positive light cone  by 
$$
\Pp = \{a\in H^2(X,\R): a^2 > 0, a\cdot[\om] > 0\},
$$
and let $\op$ be its closure.  Recall that for any
 two nonzero elements $a, a'\in
\op$, we have $a\cdot a' \ge 0$ with equality
 only if $a = \la a'$ and $a^2 = 0$.  (This
is known as the light cone lemma.)
Given $a\in H^2(X,\Z)$, we will write $\Gr(a)$ instead of $\Gr(\PD(a))$.

\begin{lemma}\label{le}  Suppose that $X$ is a symplectic manifold 
of nonsimple type.
Then for every rational class $a\in \Pp$, 
$
\Gr(qa ) \ne 0
$
for sufficiently large $q$.  Moreover if
$$
a(E) \ge 0,\;\mbox{ for all } E\in \Ee,
$$
the representation of $\PD(qa)$ by a
$J$-holomorphic curve is   connected and of multiplicity $1$.
\end{lemma}
\proof{}  When $b_2^+ = 1$ the basic fact of Taubes--Seiberg--Witten theory is
that for every class $a\in H^2(X, \Z)$ there is a number $w(a)$ called the
{\it wall-crossing } number of $a$ such that
$$
\Gr(a) \pm \Gr(K - a) = \pm w(a).
$$
In particular, if $w(a)\ne 0$ and if $[\om]\cdot(K - a)
 < 0$ then $\Gr(K-a)$ has
to be zero, so that $\Gr(a) \ne 0$.  When $b_1 = 0$, Kronheimer and
Mrowka~\cite{KM} showed that $w(a) = 1$ provided only that
$$
k(a) = -K\cdot a + a^2 \ge 0.
$$
Since $k(qa)$ is a quadratic function of $q$,  the first
 statement holds in this case.

To prove the general case first recall that when $b_2^+ = 1$ the subspace of
$H^2(X,\R)$ generated by products of elements of
 $H^1$ has dimension at most $1$. 
Moreover, when $X$ is of nonsimple type and $H^1 \ne 0$ the 
dimension is exactly
$1$.  Since $a^2 = 0$ for any $a\in (H^1)^2$, 
this subspace intersects $\op$.  Let
$a_0$ be any integral element of $(H^1)^2 \cap \op$.
Then the wall-crossing formula of
Li-Liu~\cite{LL1} says that for elements $a$ with $k(a)\ge 0$,
$w(a) \ne 0$ if and
only if $a_0\cdot(a - K/2) \ne 0$.   Since $a_0\cdot a \ne 0$ for $a\in \Pp$
by  the light cone lemma, $w(qa)\ne 0$ for large $q$.  
The desired result follows
readily.

To prove the second statement, observe that, because distinct classes in
$\Pp$ have nontrivial intersection,  the only way $\PD(a)$ can have  a
disconnected representative is if 
either $a^2 = 0$ or some of its components are
exceptional divisors.  But both of these possibilities are ruled out by our 
hypotheses. 
\QED

\subsection*{Proof of  Theorem~\ref{2}}
 If $X$
has simple type then there only are a finite number of classes $a$ with
$\Gr_0(a)\ne 0$.   Therefore $V$ is finitely generated,
by $a_1,\dots, a_p$ say. It is then clear
that there are smooth functions $c(t), \ka_i(t)\ge 0$ such that
$$
  [\om_0] = c(t) \left([\om_t] +  \sum_{ i = 1}^p\ka_i(t) a_i \right).
$$
Hence we can change $\om_t$ to an isotopy by making $p$ inflations
along the classes $a_i, i = 1,\dots, p$.

\subsection*{Proof of Theorem~\ref{1}}

Supose first that $[\om_0]$ is
rational.
Whatever the intersection form $Q_X$, there is
a basis of $H^2(X,\Q)$ formed by rational classes $n[\om_0], e_1,\dots, e_k$
with $e_j^2 < 0$ for all $j$.  Since $[\om_0]^2 > 0$ 
we may choose the integer $n$
so that the classes $n[\om_0] \pm e_j\in \Pp$.   Then,  Lemma~\ref{le} 
implies that   $\Gr(q(n[\om_0] \pm e_j))\ne 0$ for all $j$ and large $q$.  Now,
given an isotopy $\om_t$ with $[\om_0] = [\om_1]$,
 decompose its cohomology class  as
$$
[\om_t] = c(t) [\om_0] + \sum_j \la_j(t) e_j,
$$
where $c(t)> 0$.  By the openness of the set of symplectic forms, we
can  perturb $\om_t$ so that the functions $\la_j$ meet zero transversely.
Taking $a= \PD(A)$ in Lemma~\ref{infl} to be first $
q(n[\om_0] + e_1)$ and then $q(n[\om_0] - e_1)$, we  homotop $\om_t$ so
that $\la_1(t) \ge 0$ and then so that $\la_1(t) = 0$ for all $t$.
(Since $[\om_0](E) > 0$ for all $E\in \Ee$, Lemma~\ref{le}(ii) tells us that we
may ensure
that the  manifolds representing the 
duals of $q(n[\om_0] \pm e_1)$ are connected
by taking large enough $n$.)
 These homotopies will in general increase  the
function $c(t)$, but will not affect the $\la_j$ for $j\ne 1$.  Repeating this
for
$i = 2,\dots, k$, we eventually homotop $\om_t$ to a deformation such that
$[\om_t'] = c'(t) [\om_0]$.   Dividing by $c'(t)$ we find the desired isotopy.

We will deal with non rational $[\om_0]$ by following a suggestion of Biran.  
Because $\Pp$ is open, it is easy to see that
 any class $a$ in $\Pp$ is a positive sum
of rational elements of $\Pp$, that is
$$
a =   \sum_{j= 1}^{k+1} \la_j a_j, \quad \la_j \ge 0,
$$
  where each $a_j\in \Pp\cap H^2(X,\Q)$.  (Here $k = b_2^-$ as before.)
Hence we may
achieve the inflation along $qa$ by performing $k+1$ inflations along suitable
multiples of the $a_j$.\QED

\section{Embedding balls and blowing up}

We sketch the steps in the proof of Corollary~\ref{c2}.

\NI
{\bf Step 1: Normalization.} \,\, Choose an $\om$-compatible $J$ on $X$ which
is
integrable near $k$ distinct point $x_1,\dots, x_k$, and,
 for suitably small
$\de_i> 0$, fix a holomorphic and symplectic embedding
$$
\iota:\coprod_{i = 1}^k B(2\de_i) \to X
$$
which takes the center of the $i$th ball to $x_i$ for all $i$.  Given elements
$g_j, j = 0,1$ of $\Emb(\coprod B(\la_i), X)$ we may, if the $\de_i$ are small
enough, isotop these embeddings so that they both extend $\iota$.
\SS

\NI
{\bf Step 2: Forming the blow-up.}\,\,  We define the
blow up  $\Tilde X$  to be the manifold which is obtained by
cutting out the balls $\iota(B(\de_i))$ and  identifying their boundaries to
exceptional spheres $\Si_i$ via the Hopf map.  This carries a symplectic form
that integrates on $\Si_i$ to $\pi \de_i^2$.
(For more details of this step
see~\cite{INTRO,MCDPOL}.)  
It remains to define  symplectic forms on $\Tilde X$
corresponding to the $g_j$.

 For some $\nu > 0$ choose
an extension (also called $g_j$) of $g_j $ to $\coprod B(\la_i + \nu)$.  For
$i
= 1,\dots,k$ let
$$
\phi_i: B(\la_i + \nu)\to  B(\la_i + \nu)
$$
be a radial contraction which is the identity near the boundary and takes
$ B(\la_i )$ onto $B(\de_i)$ by scalar multiplication.  Define 
the map $\Phi_j :
X\to X$ for $j = 0,1$
by setting it equal to
$$
\Phi_j = g_j\circ \phi_i\circ(g_j)^{-1}\quad\mbox{on}\quad 
g_j(B(\la_i+\nu)), 1\le
i\le  k, 
$$
and extending by the identity.  Put
$$
\si_j = \Phi_j^*(\om).
$$
Then the forms $\si_0 $ and $\si_1$ both  equal the same  multiple of the 
standard form near the balls $\iota(B(\de_i))$ and so lift
 to forms which we will call
$\tsi_0, \tsi_1$ on
the  blow-up $\Tilde X$. 
  The manifold $(\TX,\tsi_j)$ is called
the blow up of $X$ by $g_j$.  Note that the weight  of the
exceptional sphere $\Si_i$ under $\tsi_j$  is
$$
\int_{\Si_i}\tsi_j = \pi \la_i^2,\quad j = 0,1.
$$ 
\SS

\NI
{\bf Step 3: Isotopies in $\TX$.}\,\,  Observe that $g_j$ may be joined to
$$
\iota' =
\iota|_{\coprod B(\de_i)}
$$
by the family of embeddings
$$
g_j^s: \coprod_{i = 1}^k B(s\la_i) \to X.
$$
Since the above blow-up construction can be done smoothly with respect to
the parameter $s$, 
there is a deformation from $\tsi_0$ to the 
blow-up of $\iota'$ and thence back to
$\tsi_1$.  Therefore, by our main result, $\tsi_0$
is isotopic to $\tsi_1$ by some isotopy $\ttau_t$.

Let $\Si$ be the union of the $k$ exceptional spheres in $\TX$.  Both $\tsi_0$
and $\tsi_1$ are equal and nondegenerate near $\Si$ by construction.
Moreover, we can arrange that the $\ttau_t$ are also nondegenerate on
$\Si$. 
To see this, observe that the deformation  between $\tsi_0$ and $\tsi_1$
that was described above consists of forms that are nondegenerate on
$\Si$.  Therefore, when constructing the inflation we can choose the
family $J_t$ so that each component of $\Si$ is $J_t$-holomorphic.  (It is 
still possible to get sufficiently generic $J_t$'s satisfying this
restriction.) Then, by using the techniques mentioned in Remark~\ref{many}
it is easy to arrange that the $\ttau_t$ are nondegenerate on $\Si$.  
Therefore, there is a family of diffeomorphisms $\tilde 
h_t: (X, \Si)\to (X,\Si)$
such that $\tilde h_t^*(\ttau_t) = \tsi_0$ on 
some neighborhood of $\Nn_\eps$ of $\Si$.  We may also assume
that $\tilde h_0 = {\rm id}$ and that $\tilde h_1 = {\rm id}$ on $\Nn_\eps$.
 Hence
there is a family $\tilde \psi_t$ of diffeomorphisms of $\TX$ such that $$
\tilde\psi_0 = \id, \quad \tilde\psi_t = \id \;\mbox{ on }\; \Nn_\eps, \quad
\tilde\psi_1^*\tilde h_1^*(\tsi_1) = \tsi_0. $$
\SS

\NI
{\bf Step 4.}\,\,  Because the forms $\ttau_t$
are constant near $\Nn_\eps$, they are the blow-ups of corresponding forms
$\tau_t$  on $X$.  Similarly, because
the  diffeomorphisms $\tilde\psi_t$ and $\tilde h_1$ are the identity near
$\Nn_\eps$, they are the blow-up of corresponding
 diffeomorphisms $\psi_t$, $h_1$ on $X$ which are the identity near the
balls $\iota(B(\de_i))$.  Therefore, by construction of $\si_0,\si_1$ we find
that
$$
\si_0 = \Phi_0^*(\om) = \psi_1^*h_1^*\Phi_1^*(\om) = \psi_1^* h_1^* \si_1.
$$
Moreover it is easy to check that
$$
F = \Phi_1\circ h_1\circ\psi_1\circ \Phi_0^{-1}
$$
is a symplectomorphism of $X$ such that $F\circ g_0 = g_1$.

Hence it remains to show that $F$ is isotopic to the identity.
But such an isotopy $F^s$ can be constructed by doing the above construction
for each $s$ where $g_j^s, s\le 1,$ is as above.  For $s$ small enough one
finds
that $g_0^s = g_1^s$ which means that $F^s$ is the identity.
There is one point worthy of note in this last step. 
 In order to construct  $F$ we
inflated a single deformation into an isotopy, using $1$-parameter families of
submanifolds.  In fact, as pointed out in Remark~\ref{many}(ii) we just need to
do one inflation process along a family of 
orthogonally intersecting submanifolds
$Z_t$. In order to inflate an arbitrary $1$-parameter 
family of deformations one
would have to use a $2$-parameter family $Z_{st}$ 
of submanifolds, and, in general,
these might encounter singularities.  However, in the 
present situation the family of
deformations is not at all arbitrary, but consists in 
lopping off the ends of the
original deformation until one arrives 
just at the center point.  Therefore, as is not
hard to check, the corresponding family of 
inflations can be constructed using the
original $1$-parameter family  $Z_t$.
 \QED

\section{Symplectic fibrations}

Let $\om_0, \om_1$ be  $\pi$-compatible forms on the oriented fibration
$\pi:X\to B$.  Our aim is to find conditions under which
these forms are either isotopic or deformation equivalent. 
Throughout this section we assume that $B$ has
dimension $2$.   At each point $p\in X$ we define $H_p $ to
be the $\om_0$-horizontal space, ie the $\om_0$-orthogonal complement to
the tangent space to the fiber through $p$.  This space has a natural
orientation. 

\begin{lemma}\label{le:1}  In the above situation, assume
 that the restriction of $\om_1$ to $H_p$ is a
nondegenerate and positive form for all $p\in X$.  Suppose further
that either $F$ has dimension $2$ or that 
the restrictions to each fiber  of  $\om_0$ and $\om_1$  are equal.
Then
the forms 
$$ 
\om_t = (1-t)\om_0 + t \om_1,\quad 0\le t \le 1,
$$
 are all
nondegenerate.  Hence if in addition  $\om_0,\om_1$
are cohomologous, they  are isotopic.
\end{lemma} 
\proof{}  This follows from a simple calculation done at
each point $p$.  Let us suppose that $\om_0$ and $\om_1$ restrict to the
form $\rho$ on the fiber $F_p$ through $p$.   Then we may choose
coordinates near $p$ so that $$ \om_0 = \rho + dx\wedge dy,\qquad
\om_1 = \rho + a dx\wedge dy + dx\wedge\al - dy\wedge \be, $$
where  $\al, \be$ are $1$-forms on $F_p$.  By hypothesis, $a >
0$.  It is easy to check that
$$
\om_1^n =   n(a+\ka)\rho^{n-1}\wedge dx
\wedge dy 
$$
where the function $\ka$ is defined by
$
(n-1)\rho^{n-2}
\al\wedge\be = \ka \rho^{n-1}.
$
Thus we must have $a + \ka > 0$.
Further,
$$
\om_t^n = n(1 -t +at  + t^2 \ka)\rho^{n-1}\wedge dx
\wedge dy,
$$
which is always nondegenerate because $|\ka|t^2 < at \le 1-t+at $ when $0\le
t\le 1$.  The calculation when $F$ has dimension $2$ is even easier.
\QED 

\begin{remark}\rm  The above argument does not always go through when
$a < 0$. For example,  $1 -t +at  + t^2 \ka < 0$ if $ a  = -7, \ka = 8, t = 1/2$.
 So the set  of cohomologous fibered forms 
need  not be convex.
\end{remark}

\begin{cor}\label{le:2} Given any  $\pi$-compatible
form $\om$, and any nondegenerate $2$-form $\si$ on the base $B$, the
forms  $\om_\mu = \om + \mu\pi^*\si$ are nondegenerate for all $\mu >
0$.  Moreover  if $H_p$ is defined with respect to any $\pi$-compatible
form $\om_0$ the restriction of $\om_\mu$ to $H_p$ is 
nondegenerate and positive when  $\mu$ is sufficiently
large.  \end{cor}  
\proof{}  This is immediate since in the previous lemma we may
identify $dx\wedge dy$ with $\pi^*(\si)$.   \QED

The above results immediately imply:

\begin{prop}\label{fibdef}  Let $\pi:X\to B$ be an oriented  fibration
with base $B$ of dimension $2$, and let  $\om_0$ and $\om_1$ be any
two $\pi$-compatible forms on $X$.  If  either $F$ has dimension $2$ or
$\om_0$ and $\om_1$ have equal restrictions to all fibers, the forms
$\om_0$ and $\om_1$ are deformation equivalent through
$\pi$-compatible forms.
\end{prop}

\subsubsection*{Proof of Proposition~\ref{fib}}

We are given a $\pi$-compatible form $\om$ om $X = F\times B$ in the
class of $ \si_F + \la\si_B$.  By Proposition~\ref{fibdef} $\om$ is
deformation equivalent to the split form $ \si_F + \la\si_B$.  Thus we
just have to convert this deformation into an isotopy.
We may assume that $B$ is not a sphere since that case is already
known.  If $B$ is also not a torus, $X$ has the smooth structure of a minimal
K\"ahler surface of general type, and so, by~\cite{TAU}
 its only non-zero Gromov invariant is
that of the canonical class $K$.  Thus $V$ is generated by the
 class
$$
\PD(K) = -c_1(TM, J) =  (2g_F - 2) [\si_F] + (2g_B - 2) [\si_B].
$$
Since our assumptions imply that
$
[\om] \in \R^+ [\om_t] + \R^+ \PD(K),
$
 the result follows from Theorem~\ref{2}.
It remains to consider the case when
 $B$ is a torus.   Taubes showed in~\cite{TAUTOR} that $\Gr(B)=2-2g_F$, so
that if $g_F > 1$,  $ V$ is generated by $[\si_F]$. Hence the
result follows as before. \QED

The above argument applies to more general fibrations.
However, it gives no
information about the manifold $T\,^4$, and also does not deal with the
question as to which symplectic forms  are isotopic to
 $\pi$-compatible forms.
\SS

 \end{document}